\RequirePackage{plautopatch}
\documentclass[12pt]{amsart}
\usepackage{lscape}
\usepackage{amssymb, enumerate,mdwlist,amsthm,amsmath,bm,ascmac,graphicx,scalerel}
\usepackage{xcolor,hyperref}

\hypersetup{
    colorlinks=true,
    citecolor=blue,
    linkcolor=red,
    urlcolor=orange,
}

\newtheorem{thm}{Theorem}
\newtheorem{prop}{Proposition}
\newtheorem{conj}{Conjecture}
\newtheorem{lemma}{Lemma}
\newtheorem{fact}{Fact}

\def\tp{\mathsf{T}}
\def\FF{\mathcal F}
\def\GG{\mathcal G}
\def\R{\mathbb R}

\def\vv{\bm{v}}
\def\one{\bm{1}}
\def\<{\langle}
\def\>{\rangle}

\usepackage{tikz}
\usetikzlibrary{positioning}
\usetikzlibrary{svg.path}
\definecolor{orcid_color}{HTML}{A6CE39}

\hypersetup{pdfborder={0 0 0}} % Remove the border of a hyperlink
\DeclareRobustCommand{\orcidicon}{%
	\raisebox{.2mm}{\scalerel*{%
	\begin{tikzpicture}[xscale=1,yscale=-1,transform shape]
	\filldraw[color=orcid_color] svg {M256,128c0,70.7-57.3,128-128,128C57.3,256,0,198.7,0,128C0,57.3,57.3,0,128,0C198.7,0,256,57.3,256,128z};
	\filldraw[color=white] svg {M86.3,186.2H70.9V79.1h15.4v48.4V186.2z} svg {M108.9,79.1h41.6c39.6,0,57,28.3,57,53.6c0,27.5-21.5,53.6-56.8,53.6h-41.8V79.1z M124.3,172.4h24.5
		c34.9,0,42.9-26.5,42.9-39.7c0-21.5-13.7-39.7-43.7-39.7h-23.7V172.4z} svg {M88.7,56.8c0,5.5-4.5,10.1-10.1,10.1c-5.6,0-10.1-4.6-10.1-10.1c0-5.6,4.5-10.1,10.1-10.1
		C84.2,46.7,88.7,51.3,88.7,56.8z};
	\end{tikzpicture}}{|}}%
}
\newcommand{\orcid}[1]{\href{https://orcid.org/#1}{\orcidicon}}

\newcommand{\arxiv}[1]{arXiv:\href{https://doi.org/10.48550/arXiv.#1}{#1}}

\title[A semidefinite programming approach to cross $2$-intersecting families]{A semidefinite programming approach to \\
cross $2$-intersecting families}

\author[Hajime Tanaka]{Hajime Tanaka\,\orcid{0000-0002-5958-0375}}
\address{\href{https://www.math.is.tohoku.ac.jp/index.html}{Research Center for Pure and Applied Mathematics}, Graduate School of Information Sciences, Tohoku University, Sendai 980-8579, Japan}
\email{htanaka@tohoku.ac.jp}
\urladdr{https://hajimetanaka.org/}
\author[Norihide Tokushige]{Norihide Tokushige\,\orcid{0000-0002-9487-7545}}
\address{College of Education, University of the Ryukyus, Nishihara  903-0213, Japan}
\email{hide@cs.u-ryukyu.ac.jp}
\urladdr{http://www.cc.u-ryukyu.ac.jp/~hide/}

\begin{document}
\maketitle
\begin{abstract}
Let $k\geq 2$ and $n\geq 3(k-1)$. Let $\mathcal{F}$ and $\mathcal{G}$ be families of $k$-element subsets of an $n$-element set. Suppose that $|F\cap G|\geq 2$ for all $F\in\mathcal{F}$ and $G\in\mathcal{G}$. We show that $|\mathcal{F}||\mathcal{G}|\leq\binom{n-2}{k-2}^2$, and determine the extremal configurations. This settles the last unsolved case of a recent result by Zhang and Wu (J.\ Combin.\ Theory Ser.\ B, 2025). We also obtain the corresponding result in the product measure setting. Our proof is done by solving semidefinite programming problems.
\end{abstract}

\section{Introduction}
Let $n$ and $t$ be positive integers, and let $[n]=\{1,2,\ldots,n\}$.
Let $2^{[n]}$ and $\binom{[n]}k$ denote the power set of $[n]$ and 
the set of $k$-element subsets of $[n]$, respectively.
We say that a family $\FF\subset 2^{[n]}$ is $t$-intersecting
if $|F\cap F'|\geq t$ for all $F,F'\in\FF$.
Erd\H{o}s, Ko, and Rado \cite{EKR} proved that a $t$-intersecting family 
$\FF\subset\binom{[n]}k$ has size at most $\binom{n-t}{k-t}$ provided $n\geq n_0(k,t)$. 
They also showed that $n_0(k,1)=2k$ is the sharp bound for $t=1$. However, finding the 
exact $n_0(k,t)$ for $t\geq 2$ is more difficult. This was done by Frankl 
\cite{Fckt} for $t\geq 15$ by 
combinatorial methods. Then Wilson \cite{W}
solved this problem for all $t$ by algebraic methods.
\begin{thm}[\cite{EKR,Fckt,W}]\label{thm:t-EKR}
Let $n\geq k\geq t\geq 1$ be integers, and let $n\geq (t+1)(k-t+1)$.
Suppose that $\FF\subset\binom{[n]}k$ is a $t$-intersecting family. Then
$|\FF|\leq\binom{n-t}{k-t}$. Moreover, if $n>(t+1)(k-t+1)$ and $|\FF|=\binom{n-t}{k-t}$,
then there exists a $t$-element subset $T\subset[n]$ such that 
$\FF=\{F\in\binom{[n]}k:T\subset F\}$.
\end{thm}

\noindent
It is much more difficult to determine the maximum size of $t$-intersecting families
for the case $n<(t+1)(k-t+1)$, and this was one of the central problems in extremal set
theory. Ahlswede and Khachatrian \cite{AK1997,AK1999} settled this problem 
(including Wilson's result) completely by combinatorial methods, and the result is
called the complete intersection theorem.

In this paper, we extend Wilson's idea to apply it for cross $t$-intersecting families.
We say that two families $\FF,\GG\subset 2^{[n]}$ are cross $t$-intersecting
if $|F\cap G|\geq t$ for all $F\in\FF$ and $G\in\GG$.
It is natural to ask the maximum of $|\FF||\GG|$ for which $\FF,\GG\subset\binom{[n]}k$
are cross $t$-intersecting families. Let $f(n,k,t)$ denote this maximum.
Pyber \cite{P} proved $f(n,k,1)=\binom{n-1}{k-1}^2$ for $n\geq 2k$ (see also \cite{MT}).
In \cite{T2013}, the following conjecture was posed.
\begin{conj}[\cite{T2013}]\label{conj-uniform}
Let $n\geq k\geq t\geq 1$ be integers, and let $n\geq (t+1)(k-t+1)$.
Suppose that $\FF,\GG\subset\binom{[n]}k$ are cross $t$-intersecting families. Then
$|\FF||\GG|\leq\binom{n-t}{k-t}^2$. Moreover, if $n>(t+1)(k-t+1)$ and 
$|\FF||\GG|=\binom{n-t}{k-t}^2$,
then there exists a $t$-element subset $T\subset[n]$ such that $\FF=\GG=\{F\in\binom{[n]}k:T\subset F\}$.
\end{conj}
\noindent
Indeed, a cheap extension of Wilson's method with Cauchy--Schwarz inequality gives us 
that $f(n,k,t)=\binom{n-t}{k-t}^2$, but only for $k/n<1-2^{-1/t}$, see \cite{T2013}. 
Very recently, Zhang and Wu \cite{ZW} proved the conjecture for $t\geq 3$,
and they also determined the extremal configurations for the case 
$|\FF||\GG|=\binom{n-t}{k-t}^2$.
Their proof is based on ideas developed by Ahlswede and Khachatrian 
\cite{AK1997}, and it seems that the smaller $t$, the more difficult the proof, as in \cite{Fckt}.
Our first result settles the last remaining case of $t=2$.
\begin{thm}\label{thm:c2k}
Let $k\geq 2$ and $n\geq 3(k-1)$.
Suppose that $\FF,\GG\subset\binom{[n]}k$ are cross $2$-intersecting families. Then
$|\FF||\GG|\leq\binom{n-2}{k-2}^2$. Moreover, if $|\FF||\GG|=\binom{n-2}{k-2}^2$, then
the following hold:
\begin{enumerate}
\item If $n>3(k-1)$, then there exists a $2$-element subset $T\subset[n]$ such that 
$\FF=\GG=\FF_0(T)$, where $\FF_0(T):=\{F\in\binom{[n]}k:T\subset F\}$.
\item If $n=3(k-1)$, then either
there exists a $2$-element subset $T\subset[n]$ such that $\FF=\GG=\FF_0(T)$,
or there exists a $4$-element subset $T'\subset[n]$
such that $\FF=\GG=\{F\in\binom{[n]}k:|F\cap T'|\geq 3\}$.
\end{enumerate}
\end{thm}
Our proof is algebraic, and it is based on ideas from \cite{ST2014BLMS}, where Suda and 
Tanaka obtained an Erd\H{o}s--Ko--Rado type result for cross $1$-intersecting families
of subspaces in a vector space. 
For the proof, they solved the corresponding semidefinite programming problem.
This SDP approach is a strong tool for intersection problems and has led to some 
successful results in \cite{ST2014BLMS, STT, TT2024}, but all of these
results are for cross $1$-intersecting families.
See also \cite{LPV2025} for some recent results on semidefinite approximations for related combinatorial problems. 

After submitting this paper, the authors learned that Chen et al.\ 
had obtained a similar but different result in \cite{CLWZ2026}. They 
proved that if $n\geq 3.38l$, $l\geq k\geq 2$,
and $\FF\subset\binom{[n]}l$ and $\GG\subset\binom{[n]}k$ are cross
2-intersecting, then $|\FF||\GG|\leq\binom{n-2}{l-2}\binom{n-2}{k-2}$.
Their proof follows the same line as in \cite{ZW}.

For a real number $p$ with $0<p<1$ and a family $\FF\subset 2^{[n]}$,
we define the $\mu_p$ measure of $\FF\subset 2^{[n]}$ by
\[
 \mu_p(\FF):=\sum_{F\in\FF}p^{|F|}(1-p)^{n-|F|}.
\]
The complete intersection theorem has its measure version, 
see \cite{DS,Friedgut2008, T2005},
which was essentially obtained by Ahlswede and Khachatrian in \cite{AK1998},
and completed by Filmus in \cite{Filmus2017}.
Here we state a special case for $p\leq\frac1{t+1}$, which can be viewed as the 
measure version of Theorem~\ref{thm:t-EKR}.
\begin{thm}[\cite{AK1998,DS,Filmus2017,Friedgut2008,T2005}]\label{thm:p-EKR}
Let $n\geq t\geq 1$ be integers, and let $p\leq\frac1{t+1}$.
Suppose that $\FF\subset 2^{[n]}$ is a $t$-intersecting family. Then
$\mu_p(\FF)\leq p^t$. Moreover, if $p<\frac1{t+1}$ and $\mu_p(\FF)=p^t$,
then there exists a $t$-element subset $T\subset[n]$ such that $\FF=\{F\in 2^{[n]}:T\subset F\}$.
\end{thm}
\noindent
Friedgut \cite{Friedgut2008} found an algebraic proof for Theorem~\ref{thm:p-EKR},
that is, a measure counterpart of Wilson's proof for Theorem~\ref{thm:t-EKR}.

In \cite{T2013CPC}, the following conjecture was posed.
\begin{conj}\label{conj-measure}
Let $n\geq t\geq 1$ be integers, and $p_1,p_2<\frac1{t+1}$.
Suppose that $\FF,\GG\subset 2^{[n]}$ are cross $t$-intersecting families. Then
$\mu_{p_1}(\FF)\mu_{p_2}(\GG)\leq(p_1p_2)^t$. Moreover, if 
$\mu_{p_1}(\FF)\mu_{p_2}(\GG)=(p_1p_2)^t$, then there exists a 
$t$-element subset $T\subset[n]$ such that $\FF=\GG=\{F\in 2^{[n]}:T\subset F\}$.
\end{conj}
\noindent
This conjecture is still wide open, and known to be true only if $t=1$ in \cite{T2010},
$p_1=p_2$ and $t\geq 14$ in \cite{FLST}, $p_1=p_2$, $t\geq 2$, and $p_1=p_2<1-2^{-1/t}$
in \cite{Filmus2013, T2013CPC}. 
(Similarly, one can consider a problem obtained by replacing the condition 
$\FF,\GG\subset\binom{[n]}k$ in Conjecture~\ref{conj-uniform} with
the condition $\FF\subset\binom{[n]}k$, $\GG\subset\binom{[n]}l$, and $k\geq l$,
and this problem is also wide open.)
Our second result verifies Conjecture~\ref{conj-measure} for the case $p_1=p_2$ and $t=2$
by the SDP approach.

\begin{thm}\label{thm:c2p}
Let $0<p\leq\frac13$. Suppose that two families $\FF,\GG\subset 2^{[n]}$
are cross $2$-intersecting. Then $\mu_p(\FF) \mu_p(\GG) \leq p^{4}$.
Moreover, if $\mu_p(\FF) \mu_p(\GG) = p^{4}$, then the following hold:
\begin{enumerate}
\item If $p<1/3$, then there exists a $2$-element subset $T\subset[n]$ such that 
$\FF=\GG=\FF_0(T)$, where $\FF_0(T):=\{F\in 2^{[n]}:T\subset F\}$.
\item If $p=1/3$, then either
there exists a $2$-element subset $T\subset[n]$ such that $\FF=\GG=\FF_0(T)$,
or there exists a $4$-element subset $T'\subset[n]$
such that $\FF=\GG=\{F\in 2^{[n]}:|F\cap T'|\geq 3\}$.
\end{enumerate}
\end{thm}

A few results are obtained in \cite{LST2017, LST2019} for the case
$p_1=p_2>\frac1{t+1}$ and large $t$.

We prove Theorem~\ref{thm:c2k} in Section~\ref{sec2}, and Theorem~\ref{thm:c2p} in 
Section~\ref{sec3}. These proofs are based on the same idea, and we solve similar
semidefinite programming problems. 
However, the actual construction of each optimal solution is quite different, and we describe them separately. In the last section, we include an
alternative proof suggested by one of the referees, and discuss how it relates to the SDP approach.

Finally, concerning cross $t$-intersecting problems,
we compare the previously known algebraic approach based on the
Hoffman's bound with Cauchy--Schwarz inequality \cite{Filmus2013, T2013, T2013CPC} and our SDP approach. 
As we mentioned, the former needs a stronger assumption
($k/n<1-2^{-1/t}$ or $p<1-2^{-1/t}$). The main reason is that this
approach uses information only between two families,
which corresponds to off-diagonal blocks in the SDP setting as we will explain
in the next section.
On the other hand, the latter SDP approach also uses information inside
each family, which corresponds to diagonal blocks.

\section{The $k$-uniform setting}\label{sec2}
\subsection{Corresponding SDP problem}\label{sec2.1}
We define a bipartite graph $G$ that represents the $t$-intersecting property.
The vertex set is $\bm\Omega=\Omega_1\sqcup\Omega_2$, where each $\Omega_i$ is a copy of 
$\binom{[n]}k$. Two vertices $x,y\in \bm\Omega$ are adjacent if and only if they come from 
different $\Omega_i$ and $|x\cap y|<t$. In this case, we write $x\sim y$. We say that 
two subsets $U_1\subset\Omega_1, U_2\subset\Omega_2$ are cross-independent if there are no 
edges between $U_1$ and $U_2$. 
This is equivalent to the condition that $U_1$ and $U_2$ are cross $t$-intersecting
families.
Our aim is to find a pair of cross-independent sets $U_1, U_2$ such that 
$|U_1||U_2|$ is maximized.

Let $\R^{\bm\Omega\times \bm\Omega}$ denote the set of real matrices indexed by 
$\bm\Omega$, and let $\R^{\bm\Omega}$ denote the set of 
real column vectors indexed by $\bm\Omega$.
For $i,j\in[2]$, 
define $\R^{\Omega_i\times\Omega_j}$ and $\R^{\Omega_i}$ similarly.
Let $J\in\R^{\Omega_i\times\Omega_j}$ be the all-ones matrix, and let
$E_{x,y}\in\R^{\Omega_i\times\Omega_j}$, where $x\in\Omega_i,y\in\Omega_j$, be the matrix whose $(x,y)$-entry is $1$,  and all the other entries are $0$. 
Let $I\in \R^{\Omega_i\times\Omega_i}$ be the identity matrix.
For $A,B\in\R^{\bm\Omega\times \bm\Omega}$, let $A\bullet B:=\operatorname{trace}(A^{\mathsf{T}}B)$ denote their trace inner product.
Finally, let $S\R^{\bm\Omega\times \bm\Omega}$ denote the set of real symmetric
matrices in $\R^{\bm\Omega\times\bm\Omega}$.

Consider the following semidefinite programming problem\footnote{The constraint $\frac12\big[\begin{smallmatrix} I & 0 \\ 0 & I \end{smallmatrix}\big] \bullet X = 1$ is weaker than $\big[\begin{smallmatrix} I & 0 \\ 0 & 0 \end{smallmatrix}\big] \bullet X = \big[\begin{smallmatrix} 0 & 0 \\ 0 & I \end{smallmatrix}\big] \bullet X= 1$ considered in \cite{ST2014BLMS} (see also \cite{STT}), but using the latter constraint instead does not decrease the optimal value. See \cite[Remark 2]{LPV2025}.} with a variable
$X \in S\mathbb{R}^{\bm\Omega \times \bm\Omega}$.
\begin{equation*}
\begin{array}{lll}
\text{(P):} & \text{maximize} & 
\dfrac12\begin{bmatrix} 0 & J\\ 
J & 0 \end{bmatrix} \bullet X \\[4mm]
& \text{subject to} & 
\dfrac12\begin{bmatrix} I & 0 \\ 0 & I \end{bmatrix} \bullet X = 1,\\[4mm]
&& \begin{bmatrix} 0 & E_{x,y} \\ 
E_{y,x} & 0 \end{bmatrix} \bullet X = 0  
\text{ for } x\in\Omega_1,\, y\in\Omega_2,\, x\sim y,\\[4mm]
&& X \succeq 0, \ X \geq 0.
\end{array}
\end{equation*}
Here, $X \succeq 0$ means that $X$ is positive semidefinite, and $X \geq 0$ means
that all the entries of $X$ are non-negative.
A feasible solution to (P) is naturally associated with every pair of cross-independent 
sets $U_1\subset\Omega_1, U_2\subset\Omega_2$ in the bipartite graph $G$.
Indeed, letting $\bm{x}_i\in\R^{\Omega_i}$ be the characteristic vector of $U_i$,
we have a feasible solution
\begin{equation}\label{eq:def of k-X}
X_{U_1,U_2}:=\begin{bmatrix} \frac{1}{\sqrt{|U_1|}}\bm{x}_1 \\ \frac{1}{\sqrt{|U_2|}}\bm{x}_2 \end{bmatrix} \!\! \begin{bmatrix} \frac{1}{\sqrt{|U_1|}}\bm{x}_1 \\ \frac{1}{\sqrt{|U_2|}}\bm{x}_2 \end{bmatrix}^{\mathsf{T}}\in S\R^{\bm\Omega\times\bm\Omega},
\end{equation}
whose objective value is
\[
\frac12\begin{bmatrix} 0 & J \\ J & 0\end{bmatrix} \bullet X_{U_1,U_2}
=\sqrt{|U_1||U_2|}.
\]

The dual problem to (P) is described as follows.
\begin{equation}\label{prob:k-D}
\begin{array}{lll}
\text{(D):} & \text{minimize} & \alpha\\
	& \text{subject to} & S:=
\dfrac12\begin{bmatrix} \alpha I & - J \\  
- J & \alpha I
\end{bmatrix} 
+{\displaystyle\sum_{x\sim y}}\,\gamma_{x,y} 
\begin{bmatrix} 0 & E_{x,y} \\ E_{y,x} & 0 \end{bmatrix} - Z \succeq 0, 
\\
		&& Z \geq 0,
\end{array}
\end{equation}
where the variables are $\alpha,\gamma_{x,y}\in\mathbb{R}$ $(x\in\Omega_1,\, y\in\Omega_2,\, x\sim y)$, and
$Z\in S\mathbb{R}^{\bm\Omega\times\bm\Omega}$.
Indeed, for any feasible solutions $X$ and $(\alpha,\gamma_{x,y}, Z)$ to (P) and (D), respectively, we have
\begin{align*}
\alpha-\frac12
\begin{bmatrix}
0&J\\
J&0
\end{bmatrix} 
\bullet X
&=
\frac12\begin{bmatrix} \alpha I & - J \\  
- J & \alpha I
\end{bmatrix} \bullet X\\
&= S\bullet X {\displaystyle-\sum_{x\sim y}}\,\gamma_{x,y} 
\begin{bmatrix} 0 & E_{x,y} \\ E_{y,x} & 0 \end{bmatrix} \bullet X + Z \bullet X\\
&=S\bullet X+Z\bullet X\geq 0.
\end{align*}
Consequently, we have the following.
\begin{prop}\label{prop-k}
Let $U_1\subset\Omega_1, U_2\subset\Omega_2$ be cross-independent sets in the bipartite graph $G$ defined above, and
let $(\alpha,\gamma_{x,y},Z)$ be a feasible solution to (D).
Then, the following hold.
\begin{enumerate}
\item We have $|U_1||U_2|\leq\alpha^2$.
\item If $|U_1||U_2|=\alpha^2$, then $S\bullet X_{U_1,U_2}=Z\bullet X_{U_1,U_2}=0$.
\end{enumerate}
\end{prop}

\subsection{Preliminaries}
Let $\Omega=\binom{[n]}k$. 
For $i=0,1$, define a matrix $B_i\in\R^{\Omega\times\Omega}$ by
\[
 (B_i)_{x,y}=\begin{cases}1&\text{if }|x\cap y|=i,\\
	      0&\text{otherwise},   \end{cases}
\]
for $x,y\in \Omega$.
These matrices are in the Bose--Mesner algebra of the Johnson scheme $J(n,k)$,
and their eigenvalues are well known. 
Here, we collect some basic facts, see, e.g., \cite{W} or \cite[Chapter 6]{GM} for the details.

For $0\leq f\leq k$, let $W_{f,k}$ and $\overline W_{\!f,k}$ be 
matrices indexed by $\binom {[n]}f\times\binom {[n]}k$ and defined by
\[
(W_{f,k})_{x,y}=\begin{cases}1&\text{if }x\subset y,\\0&\text{otherwise},\end{cases}
\quad
(\overline W_{\!f,k})_{x,y}=\begin{cases}1&\text{if }x\cap y=\emptyset,\\0&\text{otherwise}, \end{cases}
\]
for $x\in\binom {[n]}f,y\in \binom {[n]}k$.
Let $D_f=W_{f,k}^\tp \overline W_{\!f,k}$. Then, $(D_f)_{x,y}=\binom{k-|x\cap y|}f$ for $x,y\in \Omega$, and
\[
 B_0=D_k,\quad B_1=D_{k-1}-kD_k.
\]
Let $U_f$ be the column space of $W_{f,k}^\tp$, e.g., $U_0$ is spanned by
the all-ones vector $\one$.
Since $W_{j,k}^\tp W_{i,j}^\tp=\binom{k-i}{j-i}W_{i,k}^\tp$ for 
$i\leq j\leq k$, we see that $U_0\subset U_1\subset\cdots\subset U_k=\R^\Omega$.
Let $V_0=U_0$, and for $f=1,2,\dots,k$, let $V_f$ be the orthogonal
complement of $U_{f-1}$ in $U_f$, that is, $U_f=V_f\oplus U_{f-1}$.
Then, we get an orthogonal decomposition 
$\R^\Omega=V_0\oplus V_1\oplus\cdots\oplus V_k$, where 
$\dim V_f=\binom nf-\binom n{f-1}$.
Let $\{\vv^f_i:1\leq i\leq \dim V_f\}$ be an orthonormal basis of $V_f$,
and let $V$ be an $\Omega\times\binom nk$ orthogonal matrix obtained by 
arranging all these column basis vectors.
Then we have the following, see \cite{W}, or Sections 6.5 and 8.5 in \cite{GM}.

\begin{fact}\label{fact:k}
\begin{enumerate}
\item Let $0\leq f\leq k$. 
Then, $V^{-1}D_fV$ is a diagonal matrix with diagonal entries 
$(-1)^j\binom{k-j}{f-j}\binom{n-f-j}{k-j}$ with multiplicities
$\dim V_j$ for $j=0,1,\ldots,k$.
\item Let $i=0,1$. 
Then, $V^{-1}B_iV$ is a diagonal matrix with diagonal entries
$\lambda_i(j)$ with multiplicities $\dim V_j$ for $j=0,1,\ldots,k$, where
\begin{align*}
 \lambda_0(j)&=(-1)^j\binom{n-k-j}{k-j},\\
 \lambda_1(j)&=(-1)^j\binom{n-k-j}{k-j}\left(\frac{(n-k-j+1)(k-j)}{n-2k+1}-k\right).
\end{align*}
\item The only non-zero eigenvalue of $J=D_0$ is $\binom nk$ with eigenspace $U_0$. 
\end{enumerate} 
\end{fact}

We also need the following result, which is a special case of the complete intersection theorem \cite{AK1997,AK1999}, to determine the extremal configurations in Theorem~\ref{thm:c2k}.

\begin{fact}\label{fact:cit}
Let $k\geq 2$ and $n=3(k-1)$.
Suppose that $\FF\subset\binom{[n]}k$ is a $2$-intersecting family with
$|\FF|=\binom{n-2}{k-2}$.
Then, there exists a $2$-element subset $T\subset[n]$ such that 
$\FF=\{F\in\binom{[n]}k:T\subset F\}$,
or there exists a $4$-element subset $T'\subset[n]$
such that $\FF=\{F\in\binom{[n]}k:|F\cap T'|\geq 3\}$.
\end{fact}

\subsection{Proof Theorem~\ref{thm:c2k}}
Since the theorem clearly holds for $k=2$, we may assume that $k\geq 3$.
For the proof, we will construct an optimal solution to the dual problem in
\eqref{prob:k-D}. 
To this end, let $\alpha=\binom{n-2}{k-2}$ to be the optimal solution, and let
\[
 S=\frac12
\left[
\begin{matrix}
 \alpha I &- J\\
 - J&\alpha I
\end{matrix}
\right] 
-
\left[
\begin{matrix}
0 & Q\\
Q^\tp& 0 
\end{matrix}
\right] 
-
\left[
\begin{matrix}
P& 0\\
0 & P
\end{matrix}
\right],
\]
where 
\[
 P=\epsilon_0 B_0+\epsilon_1 B_1,\quad
 Q=\gamma_0 B_0+\gamma_1 B_1.
\]
Under this setting, we need to choose 
$\epsilon_0,\epsilon_1,\gamma_0,\gamma_1\in\R$ so that $\epsilon_0>0,\epsilon_1>0$,
and $S\succeq 0$.

Let $\tilde V=\left[\begin{smallmatrix} V&0\\0&V\end{smallmatrix}\right]$.
Then, by Fact~\ref{fact:k}, 
we have
\[
\tilde V^{-1} S\tilde V=
 \left[\begin{array}{cccc|cccc}
  *& & & &*\\
   &u_j& & & &v_j\\
   & &\ddots & & & &\ddots\\
   & & &* & & & & *\\
\hline
  *& & & &*\\
   &v_j& & & &u_j\\
   & &\ddots & & & &\ddots\\
   & & &* & & & & *
 \end{array}\right],
\]
that is, 
$\tilde V^{-1}S\tilde{V}$ is permutation-similar to a block diagonal matrix having diagonal blocks
$S'_j=
\left[
\begin{smallmatrix}
u_j & v_j \\ v_j & u_j 
\end{smallmatrix}
\right]$ with multiplicities $\dim V_j$ for $j=0,1,\ldots,k$, where
\begin{align*}
u_j&=\frac12\binom{n-2}{k-2}-\epsilon_0\lambda_0(j)-\epsilon_1\lambda_1(j),\\
v_j&=\delta_j-\gamma_0\lambda_0(j)-\gamma_1\lambda_1(j),\\
\delta_j&=\begin{cases}-\frac12\binom nk& \text{if }j=0,\\
	   0&\text{otherwise}.
	  \end{cases}
\end{align*}
Then, $S\succeq 0$ is equivalent to $S'_j\succeq 0$ for all $j=0,1,\ldots,k$.
\begin{lemma}\label{lem:c2k}
We can choose constants 
$\epsilon_0,\epsilon_1,\lambda_0$ and $\lambda_1$ with\footnote{The condition $\epsilon_0,\epsilon_1>0$ will be used to determine the extremal structures.}
$\epsilon_0,\epsilon_1> 0$ so that $S\succeq 0$.

\end{lemma}

\begin{proof}
We will show that $S\succeq 0$ is achieved by choosing $\epsilon_1> 0$ 
sufficiently small, and setting\footnote{
These constants are chosen so that $u_0=v_0=0$ and $u_1=-v_1$.
}
\begin{align*}
\epsilon_0&
=-\frac{k^2}{n-2k+1}\epsilon_1+\frac{k(k-1)}{2n(n-1)}\binom{n}{k}\binom{n-k}k^{-1},\\
\gamma_0&=\frac{k^2}{n-2k+1}\epsilon_1+\frac12\left(\frac{k^2(n-k)}{n(n-1)}-1\right)\binom nk\binom{n-k}k^{-1},\\
\gamma_1&=-\epsilon_1-\frac{(n-k)(n-2k+1)}{2n(n-1)}\binom nk\binom{n-k}k^{-1}.
\end{align*}

Note that if $\epsilon_1>0$ is sufficiently small, then $\epsilon_0>0$.
Note also that $S_j'\succeq 0$ is equivalent to $u_j\geq|v_j|$, 
because the eigenvalues of $S'_j$ are $u_j\pm v_j$. 
We have that
\[
 u_0=v_0=0,\quad u_1=-v_1=
\frac n{2(n-k)}\binom{n-2}{k-2}-\frac{\epsilon_1 n}{n-2k+1}\binom{n-k-1}{k-1},
\]
and so $S'_0\succeq 0$ and $S'_1\succeq 0$ (by choosing $\epsilon_1>0$ sufficiently small). 
We need to show that $S'_j\succeq 0$ for $j=2,3,\ldots, k$.

For $j\geq 2$, we have
\begin{align*}
u_j&=
\frac12\binom{n-2}{k-2}-(-1)^j C_j\left(
-\frac{j(n-j+1)}{n-2k+1}\epsilon_1+\frac{k(k-1)}{2n(n-1)}D\right),\\
v_j&=-(-1)^j C_j\left(\frac{j(n-j+1)}{n-2k+1}\epsilon_1+\frac{f(j)}{2n(n-1)}D\right),
\end{align*}
where $C_j=\binom{n-k-j}{k-j}$, $D=\binom nk\binom{n-k}k^{-1}$, and 
\begin{align*}
f(j)&=-(n-k)j^2+(n+1)(n-k)j-n(n-1).
\end{align*}
Since $\frac{\partial f(j)}{\partial j}=(n-k)(n-2j+1)>0$, 
$f(j)$ is increasing in $j$. This, together with $f(2)=(n-1)(n-2k)\geq 0$, gives us
$f(j)\geq 0$ for $2\leq j\leq k$.
Thus, by choosing $\epsilon_1>0$ sufficiently small, 
we have\footnote{
Indeed, if $\epsilon_1=0$, then 
$u_j=\frac12\binom{n-2}{k-2}(1-(-1)^j\binom{n-k-j}{k-j}\binom{n-k}k^{-1})$,
and the RHS is minimized at $j=2$, which implies that
$u_j\geq\frac12\binom{n-2}{k-2}(1-\frac{k(k-1)}{(n-k)(n-k-1)})\geq 0$.
The other estimations are obtained similarly.
} $u_j\geq 0$ for $j\geq 2$,
and $v_j \leq 0$ for even $j\geq 2$, and $v_j>0$ for odd $j\geq 3$.
Since the eigenvalues of $S_j'$ are $u_j\pm v_j$,
for proving $S_j'\succeq 0$, it suffices to show
$u_j+v_j\geq 0$ for even $j$, and $u_j-v_j\geq 0$ for odd $j$.

Suppose that $j\geq 2$ is even.
We need to show that $u_j+v_j\geq 0$, which is rewritten as
\begin{align}\label{eq:g(j)}
1\geq\binom{n-k-j}{k-j}\binom{n-k}k^{-1}\left(1+\frac{f(j)}{k(k-1)}\right)=:g(j). 
\end{align}
Indeed, we have $g(2)=1$. We also have
\[
 g(4)=\frac{(3n-k-11)(k-2)(k-3)}{(n-k-3)(n-k-2)(n-k-1)}
\]
and $(n-k-3)(n-k-2)(n-k-1)-(3n-k-11)(k-2)(k-3)=(n-3k+3)(n-4)(n-5)\geq 0$, 
from which we get $g(4)\leq 1$.
Now, we claim that $g(j)\leq 1$ not only for even integers but for all integers with 
$4\leq j\leq k$ (for later use). 
To this end, it suffices to show that $g(j)/g(j+1)\geq 1$ for $j\geq 4$,
which is equivalent to $h(j)\geq 0$ for $j\geq 4$, where\footnote{
We have $g(j)/g(j+1)\geq 1$ iff $\frac{n-k-j}{k-j}\frac{k(k-1)+f(j)}{k(k-1)+f(j+1)}\geq 1$ iff $(n-k)h(j)\geq 0$.} 
\[
h(j):= j^2 (2 k-n-2)+j n (-2 k+n+2)+2 k^2-2 k-n^2+n.
\]
Since $\frac{\partial h(j)}{\partial j}=(n-2j)(n-2k+2)>0$, 
$h(j)$ is increasing in $j$, and
\begin{align*}
h(4)&=3 n^2 -(8 k+7) n+2 k^2+30 k-32.
\end{align*}
The right hand side is minimized at $n=3(k-1)$, for which
$h(4)=(k-1)(5k-16)>0$ for $k\geq 4$. Thus, we have
$g(j)\leq 1$ for all integers $2\leq j\leq k$ except $j=3$,
in particular, $S'_j\succeq 0$ for all even $j$ with $2\leq j\leq k$.

Finally, suppose that $j\geq 3$ is odd.
We need to show that $u_j-v_j\geq 0$, which is rewritten as
\begin{align}\label{eq:g odd}
1\geq\binom{n-k-j}{k-j}\binom{n-k}k^{-1}\left(E_j\epsilon_1-1+\frac{f(j)}{k(k-1)}\right),
\end{align}
where $E_j=\frac{4n(n-1)j(n-j+1)}{k(k-1)(n-2k+1)D}$.
Comparing \eqref{eq:g odd} with \eqref{eq:g(j)}, we see that the right hand side of
\eqref{eq:g odd} is smaller than $g(j)$ provided that $\epsilon_1>0$ is sufficiently small,
and thus \eqref{eq:g odd} holds for all $j\geq 4$.

The only remaining case is $j=3$. We need to show that
\[
 1\geq\frac{k(k-1)(k-2)}{(n-k)(n-k-1)(n-k-2)}E_3\epsilon_1+
\frac{(k-2)(2n^2-(3k+5)n-k^2+7k)}{(n-k)(n-k-1)(n-k-2)}.
\]
By choosing $\epsilon_1>0$ so that $k(k-1)E_3\epsilon_1<\frac1k$, it suffices to show that
\[
1\geq
\frac{(k-2)(2n^2-(3k+5)n-k^2+7k+\frac1k)}{(n-k)(n-k-1)(n-k-2)}, 
\]
that is,
\[
F(n):=n^3-(5k-1) n^2+\left(6 k^2+5 k-8\right) n-12k^2+12k+\frac{2}{k}-1\geq 0
\]
Since $\frac{\partial F(n)}{\partial n}>0$ for $n\geq 3(k-1)$ and $k\geq 3$, 
it suffices to check it at $n=3(k-1)$, and 
\[
 F(3(k-1))=\frac1k(k-2)(3k^2-3k-1)>0
\]
for $k\geq 3$. 
This completes the proof of the lemma.
\end{proof}

By Lemma~\ref{lem:c2k} we have constructed an optimal solution to \eqref{prob:k-D}, and the inequality of the theorem follows from Proposition~\ref{prop-k} (i). 

To determine the extremal configurations, 
let $n\geq 3(k-1)$, and suppose that $U_1,U_2\subset \binom{[n]}k$ 
are cross $2$-intersecting families with $|U_1||U_2|=\binom{n-2}{k-2}^2$. 
Let $X=X_{U_1,U_2}$ be defined by \eqref{eq:def of k-X}. Then, by Proposition~\ref{prop-k} (ii), 
we have $Z\bullet X=0$, where
\[
 Z=\begin{bmatrix}P&0\\0&P\end{bmatrix}
=\begin{bmatrix}\epsilon_0 B_0+\epsilon_1 B_1&0\\
0&\epsilon_0 B_0+\epsilon_1 B_1\end{bmatrix}.
\]
By our construction, it follows that 
$(\epsilon_0 B_0+\epsilon_1 B_1)_{x,y}>0$ if $|x\cap y|\leq 1$.
Then, $Z\bullet X=0$ implies that $(X)_{x,y}=0$ if $|x\cap y|\leq 1$, that
is, each $U_i$ ($i=1,2$) is a $2$-intersecting family by itself. 
Thus, the corresponding configurations are determined by 
Theorem~\ref{thm:t-EKR} and Fact~\ref{fact:cit}.
\qed

\section{The measure setting}\label{sec3}
\subsection{Corresponding SDP problem}
We define a bipartite graph $G$ similarly as in Section~\ref{sec2.1}.
The only difference is that $\Omega_i$ is a copy of $2^{[n]}$ instead of $\binom{[n]}k$.
We will find a pair of cross independent sets $U_1\subset\Omega_1$ and 
$U_2\subset\Omega_2$ such that $\mu_p(U_1)\mu_p(U_2)$ is maximized.
Let $\Delta\in \R^{\Omega_i\times\Omega_i}$ be the diagonal matrix whose
$(x,x)$-entry is $\mu_{p}(\{x\})$, where $x\in \Omega_i$.

Consider the following semidefinite programming problem with a variable
$X \in S\mathbb{R}^{\bm\Omega \times \bm\Omega}$.
\begin{equation*}
\begin{array}{lll}
\text{(P):} & \text{maximize} & 
\dfrac12\begin{bmatrix} 0 & \Delta J\Delta \\ 
\Delta J\Delta & 0 \end{bmatrix} \bullet X \\[4mm]
& \text{subject to} & 
\dfrac12\begin{bmatrix} \Delta & 0 \\ 0 & \Delta \end{bmatrix} \bullet X = 1,\\[4mm]
&& \begin{bmatrix} 0 & E_{x,y} \\ 
E_{y,x} & 0 \end{bmatrix} \bullet X = 0  
\text{ for } x\in\Omega_1,\, y\in\Omega_2,\, x\sim y,\\[4mm]
&& X \succeq 0, \ X \geq 0.
\end{array}
\end{equation*}
Again, a feasible solution to (P) is associated with every pair of cross-independent 
sets $U_1\subset\Omega_1,U_2\subset\Omega_2$ in the bipartite graph $G$.
Indeed, letting $\bm{x}_i\in\R^{\Omega_i}$ be the characteristic vector of $U_i$,
we have a feasible solution
\begin{equation}\label{eq:def of X}
X_{U_1,U_2}:=\begin{bmatrix} \frac{1}{\sqrt{\mu_{p}(U_1)}}\bm{x}_1 \\ \frac{1}{\sqrt{\mu_{p}(U_2)}}\bm{x}_2 \end{bmatrix} \!\! \begin{bmatrix} \frac{1}{\sqrt{\mu_{p}(U_1)}}\bm{x}_1 \\ \frac{1}{\sqrt{\mu_{p}(U_2)}}\bm{x}_2 \end{bmatrix}^{\mathsf{T}}\in S\R^{\bm\Omega\times\bm\Omega},
\end{equation}
whose objective value is
\[
\frac12\begin{bmatrix} 0 & \Delta J\Delta \\ \Delta J\Delta & 0\end{bmatrix} \bullet X_{U_1,U_2}
=\sqrt{\mu_{p}(U_1)\mu_{p}(U_2)}.
\]

The dual problem to (P) is described as follows.
\begin{equation}\label{prob:D}
\begin{array}{lll}
\text{(D):} & \text{minimize} & \alpha\\
	& \text{subject to} & S:=
\dfrac12\begin{bmatrix} \alpha \Delta & - \Delta J \Delta \\  
- \Delta J\Delta & \alpha \Delta 
\end{bmatrix} 
+{\displaystyle\sum_{x\sim y}}\,\gamma_{x,y} 
\begin{bmatrix} 0 & E_{x,y} \\ E_{y,x} & 0 \end{bmatrix} - Z \succeq 0, 
\\
		&& Z \geq 0,
\end{array}
\end{equation}
where the variables are $\alpha,\gamma_{x,y}\in\mathbb{R}$ $(x\in\Omega_1,\, y\in\Omega_2,\, x\sim y)$, and
$Z\in S\mathbb{R}^{\bm\Omega\times\bm\Omega}$.
Indeed, for any feasible solutions $X$ and $(\alpha,\gamma_{x,y}, Z)$ to (P) and (D), respectively, we have
\begin{align*}
\alpha-\frac12
\begin{bmatrix}
0&\Delta J\Delta\\
\Delta J\Delta&0
\end{bmatrix} 
\bullet X
&=
\frac12\begin{bmatrix} \alpha \Delta & - \Delta J \Delta \\  
- \Delta J \Delta & \alpha \Delta
\end{bmatrix} \bullet X\\
&= S\bullet X {\displaystyle-\sum_{x\sim y}}\,\gamma_{x,y} 
\begin{bmatrix} 0 & E_{x,y} \\ E_{y,x} & 0 \end{bmatrix} \bullet X + Z \bullet X\\
&=S\bullet X+Z\bullet X\geq 0.
\end{align*}
Consequently, we have the following.
\begin{prop}\label{prop}
Let $U_1\subset\Omega_1, U_2\subset\Omega_2$ be cross-independent sets in the bipartite graph $G$ defined above, and
let $(\alpha,\gamma_{x,y},Z)$ be a feasible solution to (D).
Then, the following hold.
\begin{enumerate}
\item We have $\mu_{p}(U_1)\mu_{p}(U_2)\leq \alpha^2$.
\item If $\mu_{p}(U_1)\mu_{p}(U_2)=\alpha^2$, then $S\bullet X_{U_1,U_2}=Z\bullet X_{U_1,U_2}=0$.
\end{enumerate}
\end{prop}

\subsection{Preliminaries}
Let $p$ be a real number with $0<p\leq \frac13$, and let $q:=1-p$.
Let 
\begin{align*}
A'&=\left[
\begin{matrix}1-\frac pq &\frac pq \\1&0\end{matrix}\right],&
I'&=\left[\begin{matrix}1&0\\0&1\end{matrix}\right],&
J'&=\begin{bmatrix} 1 & 1 \\ 1 & 1 \end{bmatrix},\\
D'&=\begin{bmatrix} 1 & 0 \\ 0 & -p/q \end{bmatrix}, &
V'&=\begin{bmatrix} 1 & \sqrt{p/q} \\ 1 & -\sqrt{q/p} \end{bmatrix}, &
\Delta'&=\begin{bmatrix} q & 0 \\ 0 & p \end{bmatrix},
\end{align*}
where the rows and columns are indexed in the order $0,1$.

For matrices 
$M=\left[\begin{smallmatrix}m_{0,0}&m_{0,1}\\m_{1,0}&m_{1,1}\end{smallmatrix}\right]$ 
and $M'$, we define the Kronecker product 
$M\otimes M':=
\left[\begin{smallmatrix}m_{0,0}M'&m_{0,1}M'\\m_{1,0}M'&m_{1,1}M'\end{smallmatrix}\right]$, 
and the $n$-fold product $M^{\otimes n}:=M\otimes M\otimes\cdots\otimes M$ ($n$ copies). 
Then, we can view the rows and columns of $M^{\otimes n}$ as indexed by the $n$-bit strings,
and the value of the $(x,y)$-entry of $M^{\otimes n}$ is $\prod_{i=1}^n m_{x_i,y_i}$.

Let $D=D'^{\otimes n}$ and
$V=V'^{\otimes n}$.
We identify the elements of $\Omega=2^{[n]}$ with their characteristic vectors in $\{0,1\}^n$, and view these matrices as matrices in $\R^{\Omega\times\Omega}$.
We note that $I=I'^{\otimes n}$, $J=J'^{\otimes n}$, and $\Delta=\Delta'^{\otimes n}$.
For $i=0,1$, we also define matrices $B_i,D_i\in\R^{\Omega\times\Omega}$ as follows. 
Let $B_0=A'^{\otimes n}$, $D_0=D'^{\otimes n}$, and
\[
 B_1=\frac1n\sum_{j=1}^n\bigotimes_{i=1}^nA'_{i,j}, \quad
 D_1=\frac1n\sum_{j=1}^n\bigotimes_{i=1}^nD'_{i,j},
\]
where
\[
 A'_{i,j}=\begin{cases}I'&\text{if }i=j,\\A'&\text{otherwise},\end{cases}\quad
 D'_{i,j}=\begin{cases}I'&\text{if }i=j,\\D'&\text{otherwise}.\end{cases}
\]
For example, if $n=3$, then
$B_1=\frac13\left(
(I'\otimes A'\otimes A')+
(A'\otimes I'\otimes A')+
(A'\otimes A'\otimes I')
\right)
$. 
The matrices $B_i$ can be defined for $i\geq 2$ as well, 
but we will not use them. These matrices were introduced by Friedgut
to give an alternative proof of Theorem~\ref{thm:p-EKR}, and used by
Filmus to prove Conjecture~\ref{conj-measure} for the case $p<1-2^{-1/t}$, see
\cite{Filmus2013} for more details.
We also note that the above matrices in $\R^{\Omega\times\Omega}$ all belong to the Terwilliger algebra of the binary Hamming scheme $H(n,2)$ \cite{Go2002EJC, Terwilliger1992JAC}.

We list some basic properties concerning these matrices.
\begin{fact}\label{fact:p}
\begin{enumerate}
\item $\Delta,\Delta J\Delta,\Delta B_0,\Delta B_1\in S\R^{\Omega\times\Omega}$.
\item\label{fact:p2} $V^\tp \Delta V=I$, $V^\tp (\Delta J\Delta) V=E_{\emptyset,\emptyset}$. 
\item For $i=0,1$, $V^\tp (\Delta B_i) V=D_i$.
\item For $i=0,1$, $(D_i)_{x,x}=\lambda_i(|x|)$, where
\begin{align*}
	\lambda_0(j):=\left(-\frac{p}{q}\right)^j, \quad \lambda_1(j):=\left(-\frac{p}{q}\right)^j\left(1-\frac{j}{np}\right).
\end{align*}
\item $(B_0)_{x,y}=
\begin{cases}\left(1-\frac pq\right)^{n-|x|-|y|}\left(\frac pq\right)^{|y|}& \text{if }|x\cap y|=0,\\
0&\text{otherwise}. \end{cases}$

\item $(B_1)_{x,y}=
\begin{cases}\frac{n-|x|-|y|}n\left(1-\frac pq\right)^{n-1-|x|-|y|}\left(\frac pq\right)^{|y|}& \text{if $|x\cap y|=0$ and $|x|+|y|<n$},\\[2mm]
\frac1n\left(1-\frac pq\right)^{n+1-|x|-|y|}\left(\frac pq\right)^{|y|-1}& \text{if $|x\cap y|=1$},\\
0&\text{otherwise}. \end{cases}$

\item\label{fact:p7} If $\epsilon_0>0$ and $\epsilon_1>0$, then
$(\epsilon_0\Delta B_0+\epsilon_1\Delta B_1)_{x,y}>0$ for all 
$x,y$ with $|x\cap y|\leq 1$.
\end{enumerate}
\end{fact}

We also need the following result, which is a special case of a result due to Filmus \cite{Filmus2017}, to determine the extremal configurations in Theorem~\ref{thm:c2p}.

\begin{fact}\label{fact:filmus}
Let $p\leq 1/3$.
Suppose that $\FF\subset 2^{[n]}$ is a $2$-intersecting family with
$\mu_p(\FF)=p^2$.
Then, there exists a $2$-element subset $T\subset[n]$ such that 
$\FF=\{F\in 2^{[n]}:T\subset F\}$,
or there exists a $4$-element subset $T'\subset[n]$
such that $\FF=\{F\in 2^{[n]}:|F\cap T'|\geq 3\}$.
\end{fact}

\subsection{Proof Theorem~\ref{thm:c2p}}
For the proof, we will construct an optimal solution to the dual problem in
\eqref{prob:D}. 
To this end, let $\alpha=p^2$, and let
\[
 S=\frac12
\left[
\begin{matrix}
 \alpha \Delta &- \Delta J\Delta\\
 - \Delta J\Delta&\alpha \Delta 
\end{matrix}
\right] 
-
\left[
\begin{matrix}
0 & Q\\
Q^\tp& 0 
\end{matrix}
\right] 
-
\left[
\begin{matrix}
P& 0\\
0 & P
\end{matrix}
\right],
\]
where 
\[
 P=\epsilon_0\Delta B_0+\epsilon_1\Delta B_1,\quad
 Q=\gamma_0\Delta B_0+\gamma_1\Delta B_1.
\]
Under this setting, we need to choose 
$\epsilon_0,\epsilon_1,\gamma_0,\gamma_1\in\R$ so that $\epsilon_0>0,\epsilon_1\geq 0$,
and $S\succeq 0$. (We will choose $\epsilon_1=0$ for $p=\frac13$ and
$\epsilon_1>0$ for $p<\frac13$.)

The matrix $S$ is divided into four large blocks corresponding to
$\Omega_1\times\Omega_1$, $\Omega_1\times\Omega_2$,
$\Omega_2\times\Omega_1$, and $\Omega_2\times\Omega_2$, respectively.
Each block is diagonalized as follows.\footnote{
Indeed, letting $\tilde \Delta=\left[\begin{smallmatrix} \Delta&0\\0&\Delta\end{smallmatrix}\right]$, each block of $\tilde\Delta^{-\frac12}S\tilde\Delta^{-\frac12}$ is 
diagonalized by the orthogonal matrix $\tilde W:=\tilde V^\tp\tilde \Delta^{\frac12}$, that is, 
$\tilde W(\tilde \Delta^{-\frac12}S\tilde \Delta^{-\frac12})\tilde W^{-1}=
\tilde V^\tp S \tilde V$.}
\begin{equation}\label{block-wise diagonalization}
S':=\tilde V^\tp S\tilde V=
 \left[\begin{array}{cccc|cccc}
  *& & & &*\\
   &(S')_{x,x}& & & &(S')_{x,y}\\
   & &\ddots & & & &\ddots\\
   & & &* & & & & *\\
\hline
  *& & & &*\\
   &(S')_{y,x}& & & &(S')_{y,y}\\
   & &\ddots & & & &\ddots\\
   & & &* & & & & *
 \end{array}\right],
\end{equation}
where 
$\tilde V=\left[\begin{smallmatrix} V&0\\0&V\end{smallmatrix}\right]$.
To show that $S\succeq 0$, it suffices to show that $S'\succeq 0$. 
For $x\in\Omega_1$ and $y\in\Omega_2$, where $x=y$ as subsets of $[n]$, 
we assign a $2\times 2$ principal submatrix of $S'$, that is,
\[
 \begin{bmatrix}
  (S')_{x,x}&  (S')_{x,y}\\
  (S')_{y,x}&  (S')_{y,y}
 \end{bmatrix}.
\]
This matrix only depends on $j=|x|=|y|$, so we define a symmetric matrix
$S'_j:=
\left[
\begin{smallmatrix}
u_j & v_j \\ v_j & u_j 
\end{smallmatrix}
\right]$ by
\begin{align*}
u_j&=\frac{p^2}2-\epsilon_0\lambda_0(j)-\epsilon_1\lambda_1(j),\\
v_j&=\delta_j-\gamma_0\lambda_0(j)-\gamma_1\lambda_1(j),\\
\delta_j&=\begin{cases}-\frac12& \text{if }j=0,\\
	   0&\text{otherwise}.
	  \end{cases}
\end{align*}
Then, $S'\succeq 0$ is equivalent to $S'_j\succeq 0$ for all $j=0,1,\ldots,n$.
\begin{lemma}\label{lem:c2p}
We can choose constants 
$\epsilon_0,\epsilon_1,\lambda_0$ and $\lambda_1$ with 
$\epsilon_0>0$ and $\epsilon_1\geq 0$ so that $S\succeq 0$, where $\epsilon_1=0$ when $p=\frac{1}{3}$ and $\epsilon_1>0$ when $p<\frac{1}{3}$.
\end{lemma}

\begin{proof}
We will show that $S'\succeq 0$ is achieved by choosing $\epsilon_1\geq 0$ 
sufficiently small, and letting
\begin{align*}
\epsilon_0&=\frac{p^2}2-\epsilon_1,\\
\gamma_0&=-\frac12+\frac12pqn+\epsilon_1,\\
\gamma_1&=-\frac12pqn-\epsilon_1.
\end{align*}
In this case, by choosing $0\leq\epsilon_1<\frac{p^2}2$, it follows that $\epsilon_0>0$.

Recall that $S_j'\succeq 0$ is equivalent to $u_j\geq|v_j|$. 
We have that
\[
 u_0=v_0=0,\quad u_1=-v_1=\frac{p^2n-2\epsilon_1}{2qn}.
\]
Thus, $S_0'\succeq 0$ and $S_1'\succeq 0$ (provided $\epsilon_1\geq 0$ is sufficiently small). 
We need to show that $S_j'\succeq 0$ for $j=2,3,\ldots, n$.
For $j\geq 2$ we have
\begin{align*}
u_j&=
\frac1{2p}\left(p^3-\left(-\frac pq\right)^j\left(p^3-\frac{2\epsilon_1j}n\right)\right),\\
v_j&=-\frac1{2p}\left(-\frac pq\right)^j\left(p(qj-1)+\frac{2\epsilon_1j}n\right).
\end{align*}

First, we deal with the case $p=\frac13$. In this case we have
\begin{equation}\label{p=1/3 uj+vj}
u_j+v_j=
\frac{1}{9\cdot 2^{j+1}}  \left(-6 (-1)^j j+8 (-1)^j+2^j\right) \geq 0.
\end{equation}
For $u_j-v_j$, we note that
\[
 u_3-v_3=-\frac{9\epsilon_1}{4n}.
\]
Thus, for $u_3-v_3\geq 0$, we necessarily have that $\epsilon_1=0$, and
in this case, we indeed have that
\begin{equation}\label{p=1/3 uj-vj}
 u_j-v_j=\frac1{18}+\frac{(-1)^j}{9\cdot 2^j}(3j-5)\geq 0
\end{equation}
for all $2\leq j\leq n$. Therefore, if $p=\frac13$, then $S_j'\succeq 0$ for all $j\geq 2$
by choosing 
$\epsilon_1=0$.

Next, we deal with the case $p<\frac13$.
By choosing $\epsilon_1>0$ sufficiently small, we have $u_j>0$. We also have
$v_j<0$ if $j$ is even, and $v_j>0$ if $j$ is odd.
Thus, to show that $S'_j\succeq 0$, it suffices to verify that
$u_j+v_j\geq 0$ if $j$ is even, and $u_j-v_j\geq 0$ if $j$ is odd.

Suppose that $j\geq 2$ is even, and let 
$f(j)=(\frac pq)^j(j-1-p)$. Then, 
\[
 u_j+v_j=\frac12\left(p^2-q f(j)\right),
\]
and $f(j)$ is positive and maximized at $j=2$. Indeed, since
$\frac{f(j)}{f(j+2)}=\frac{q^2(j-1-p)}{p^2(j+1-p)}$ and
\[
q^2(j-1-p)-p^2(j+1-p)=
j(1-2p)-q\geq 2(1-2p)-q=1-3p>0, 
\]
it follow that $f(j)>f(j+2)$.
Hence, $u_j+v_j\geq\frac12\left(p^2-q f(2)\right)=0$.

Now, suppose that $j\geq 3$ is odd, and let $f(j)=(\frac pq)^j(jq-1-p^2)$. Then,
\[
 u_j-v_j=\frac1{2}\left(p^2-f(j)-\left(\frac{p}{q}\right)^j\frac{4j\epsilon_1}{np}\right).
\]
We claim that $f(j)\leq f(3)$. 
Indeed, we have $\frac{f(j)}{f(j+2)}=\frac{q^2(jq-1-p^2)}{p^2((j+2)q-1-p^2)}$, and
\begin{align*}
q^2(jq-1-p^2)-p^2((j+2)q-1-p^2)
&=j (1-p)(1-2p)+4 p^3-3 p^2+2 p-1\\
&\geq 3(1-p)(1-2p)+4 p^3-3 p^2+2 p-1 >0
\end{align*}
for $0<p<\frac13$.
Then, choosing $\epsilon_1>0$ sufficiently small, we get
\begin{align*}
p^2-f(j)\geq p^2-f(3)=p^2(1-2p)(1-3p)/q^3>\left(\frac pq\right)^3\frac{4\epsilon_1}{p}
\geq\left(\frac pq\right)^j\frac{4j\epsilon_1}{np},
\end{align*}
which means $u_j-v_j\geq 0$.
This completes the proof of the lemma.
\end{proof}

By Lemma~\ref{lem:c2p}, we have constructed an optimal solution to \eqref{prob:D}, and the inequality of the theorem follows from Proposition~\ref{prop} (i). 

It remains to describe the extremal configurations.
Suppose that $U_1,U_2\subset 2^{[n]}$ 
are cross $2$-intersecting families with $\mu_p(U_1)\mu_p(U_2)=p^4$. Let $X=X_{U_1,U_2}$ be defined by 
\eqref{eq:def of X}. Then, by Proposition~\ref{prop} (ii), we have $Z\bullet X=0$, where
\[
 Z=\begin{bmatrix}P&0\\0&P\end{bmatrix}
=\begin{bmatrix}\epsilon_0\Delta B_0+\epsilon_1\Delta B_1&0\\
0&\epsilon_0\Delta B_0+\epsilon_1\Delta B_1\end{bmatrix}.
\]

First, assume that $p<\frac13$.
Then, we have $\epsilon_0>0$ and $\epsilon_1>0$, and it follows from
Fact~\ref{fact:p}\,(vii) that each $U_i$ ($i=1,2$) is a $2$-intersecting family by itself. 
Thus, the corresponding configurations are determined by 
Theorem~\ref{thm:p-EKR} and Fact~\ref{fact:filmus}.

Next, assume that $p=\frac{1}{3}$.
We have $\epsilon_0>0$ and $\epsilon_1=0$ in this case, and this only shows that each $U_i$ ($i=1,2$) is a $1$-intersecting family.
For $i=1,2$ and $x\in 2^{[n]}=\Omega_i$, let $\omega_i(x)$ denote the $x$-entry of $\frac{1}{\sqrt{\mu_{p}(U_i)}}V^{-1}\bm{x}_i=\frac{1}{\sqrt{\mu_{p}(U_i)}}V^{\tp}\Delta\bm{x}_i$, see Fact~\ref{fact:p}\,(ii).
In particular, we have $\omega_i(\emptyset)=\sqrt{\mu_p(U_i)}$ because the $\emptyset$-column of $V$ is the all one's vector. 
Moreover, since
\begin{equation*}
	(V^{-1}\bm{x}_i)^{\tp}(V^{-1}\bm{x}_i)=\bm{x}_i^{\tp}\Delta VV^{-1}\bm{x}_i=\bm{x}_i^{\tp}\Delta\bm{x}_i=\mu_{p}(U_i),
\end{equation*}
we have $\sum_{x\in 2^{[n]}}\omega_i(x)^2=1$ for $i=1,2$.

By Proposition~\ref{prop} (ii) and \eqref{block-wise diagonalization}, we have
\begin{align*}
	0&= S\bullet X = S'\bullet (\tilde{V}^{-1}X\tilde{V}^{-\tp}) = \sum_{x\in 2^{[n]}} S_{|x|}'\bullet \begin{bmatrix} \omega_1(x) \\ \omega_2(x) \end{bmatrix}\!\! \begin{bmatrix} \omega_1(x) \\ \omega_2(x) \end{bmatrix}^{\tp}. 
\end{align*}
Since $S'_j\succeq 0$ for all $j=0,1,\dots,n$, each of the summands in the RHS above equals zero.
Recall that $u_0=v_0=0$ and $u_1=-v_1=\frac{p^2}{2q}=\frac{1}{12}$.
Moreover, for $j\geq 2$, it follows from \eqref{p=1/3 uj+vj} and \eqref{p=1/3 uj-vj} that $u_j+v_j=0$ if and only if $j\in\{2,4\}$,
and that $u_j-v_j=0$ if and only if $j=3$.
By these comments, we have $\omega_1(x)=\omega_2(x)$ when $|x|\in\{1,2,4\}$, $\omega_1(x)=-\omega_2(x)$ when $|x|=3$, and $\omega_1(x)=\omega_2(x)=0$ when $|x|\geq 5$.
In particular, this shows that
\begin{equation*}
	\mu_p(U_1)=1-\sum_{x\ne \emptyset}\omega_1(x)^2=1-\sum_{x\ne \emptyset}\omega_2(x)^2=\mu_p(U_2),
\end{equation*}
and so $\mu_p(U_1)=\mu_p(U_2)=p^2=\frac{1}{9}$.

We now claim that $\omega_i(x)=0$ $(i=1,2)$ when $|x|=3$.
On the one hand, since $U_1,U_2$ are cross $2$-intersecting (and hence cross $1$-intersecting), we have
\begin{equation*}
	0=\frac{1}{p^2} \bm{x}_1^{\tp}(\Delta B_0)\bm{x}_2 = \frac{1}{p^2}(V^{-1}\bm{x}_1)^{\tp}D_0(V^{-1}\bm{x}_2)=\sum_{x\in 2^{[n]}}\lambda_0(|x|)\omega_1(x)\omega_2(x)=\sigma-\rho
\end{equation*}
by Fact~\ref{fact:p}\,(iii), where $\sigma:=\sum_{|x|\ne 3}\lambda_0(|x|)\omega_1(x)^2$ and $\rho:=\lambda_0(3)\sum_{|x|=3}\omega_1(x)^2$.
On the other hand, since $U_1$ is $1$-intersecting, we similarly have
\begin{equation*}
	0=\frac{1}{p^2} \bm{x}_1^{\tp}(\Delta B_0)\bm{x}_1 =\sigma+\rho.
\end{equation*}
It follows that $\sigma=\rho=0$, and the claim follows since $\lambda_0(3)\ne 0$.
This establishes that $\omega_1(x)=\omega_2(x)$ for all $x\in 2^{[n]}$, and hence we have $\bm{x}_1=\bm{x}_2$, or equivalently, $U_1=U_2$.
It follows that $U_1$ and $U_2$ are identical optimal $2$-intersecting families and are determined by Fact~\ref{fact:filmus}.

This completes the proof of Theorem~\ref{thm:c2p}.
\qed

\section*{Acknowledgments}
The authors thank one of the referees for giving an alternative proof (described in terms of Fourier analysis)
and for pointing out a result in \cite{Filmus2017} that makes it possible to determine
the extremal structure in Theorem~\ref{thm:c2p}. 
They also appreciate the other referee for many suggestions to improve the presentation.

HT was supported by JSPS KAKENHI Grant Number JP23K03064.
NT was supported by JSPS KAKENHI Grant Number JP23K03201.

\end{document}